# Bayesian sample sizes for exploratory clinical trials comparing multiple experimental treatments with a control

John Whitehead[a][*][+], Faye Cleary[b] and Amanda Turner[a]
[a]*Department of Mathematics and Statistics, Lancaster University, UK*
[b]*GlaxoSmithKline, Stockley Park, UK*

*11 July 2014*

## SUMMARY

In this paper, a Bayesian approach is developed for simultaneously comparing multiple experimental treatments with a common control treatment in an exploratory clinical trial. The sample size is set to ensure that, at the end of the study, there will be at least one treatment for which the investigators have a strong belief that it is better than control, or else they have a strong belief that none of the experimental treatments are substantially better than control. This criterion bears a direct relationship with conventional frequentist power requirements, while allowing prior opinion to feature in the analysis with a consequent reduction in sample size. If it is concluded that at least one of the experimental treatments shows promise, then it is envisaged that one or more of these promising treatments will be developed further in a definitive phase III trial.

    The approach is developed in the context of normally distributed responses sharing a common standard deviation regardless of treatment. To begin with, the standard deviation will be assumed known when the sample size is calculated. The final analysis will not rely upon this assumption, although the intended properties of the design may not be achieved if the anticipated standard deviation turns out to be inappropriate. Methods that formally allow for uncertainty about the standard deviation, expressed in the form of a Bayesian prior, are then explored.

    Illustrations of the sample sizes computed from the new method are presented, and comparisons are made with frequentist methods devised for the same situation.

*\*Correspondence to: John Whitehead, Department of Mathematics and Statistics, Fylde College, Lancaster University, Lancaster LA1 4YF, UK*
[+]*E-mail [j.whitehead@lancaster.ac.uk](mailto:j.whitehead@lancaster.ac.uk)*

# 1. Introduction

A commonly occurring situation in clinical research is the need to evaluate simultaneously alternative new approaches to the treatment of the same condition. In the pharmaceutical industry, the alternatives might be different doses, formulations or administration schedules of the same experimental drug, or they might relate to the use of that drug with and without one or more concurrent medications. In the public sector, a head-to-head comparison of several drugs or other forms of treatment from different sources might be of interest. The method described here considers the treatments to be qualitatively different. For example, if they do represent different doses of the same drug, no assumptions are made concerning the dose-response relationship, and if they are different combinations of the same set of drugs, no interactions are ignored. A clinical trial is considered in which patients are randomised between k new treatments and a control treatment. The latter might comprise or include administration of a placebo, otherwise it will often represent a standard form of treatment. Often, multiple treatments are assessed during phase II of a drug development programme, with a view to selecting one or more of them for further investigation in phase III. This paper is written in that setting, although the methods presented might be applicable in different contexts. In this paper, the whole trial is completed in a single stage, although the methodology would lend itself to a multi-stage approach.

Frequentist sample size calculations for multiple treatment comparisons are usually based on considerations of type I error and power. Type I error is the probability of claiming that one or more of the experimental treatments is better than control (and therefore worthy of further investigation), when in fact all treatments have the same effect. The power will be the probability of reaching the same conclusion under a specific scenario in which one of the new treatments is indeed superior to placebo by a specified margin. An attractive frequentist solution to this problem was suggested by Dunnett [1], and this approach has been implemented in trials such as an anti-hypertensive study described by Carlsen et al. [2] (data from which are used for an illustration in Section 5 below) and generalised to multi-stage trials [3, 4]. Other frequentist sample size calculations for single-stage multiple-treatment trials have also been suggested [5, 6]. These are based on a final F-test to determine whether a difference between the treatments exists, and are unsuitable when one of the treatments is a control and it is important to know which treatments are the most efficacious.

In this paper, we develop a Bayesian approach to sample size determination for multiple treatment trials. It is an extension of the method described by Whitehead et al. [7], extending an approach mentioned by Simon [8], to comparisons of a single experimental treatment with control. Reference [7] includes a review of a range of Bayesian methods that have been suggested for setting sample size. Spiegelhalter et al. [9], classify such approaches as "hybrid classical and Bayesian", "proper Bayesian" and "decision-theoretic Bayesian". The first of these uses Bayesian methodology only for the design of the trial, assuming that the final analysis will be frequentist: an example is provided by O'Hagan and Stevens [10]. Proper Bayesian methods are usually based on the predictive power of a final Bayesian analysis. This paper adopts a slightly different approach which is nevertheless "proper Bayesian". Finally, another Bayesian option for setting sample sizes is to use decision theory [11, 12].

The method of [7] concerns a parameter $\delta$ representing the advantage of a single experimental treatment over the control. While the value $\delta = 0$ represents no difference between the treatments, the value $\delta = \delta^*$ (where $\delta^* > 0$) would reflect a clinically important treatment effect. The Bayesian criterion for the required sample size starts with a prior distribution for $\delta$, which will be transformed into a posterior distribution after observing the trial data. If the posterior belief that $\delta > 0$ is large enough, then it will be concluded that the experimental treatment is promising and research will proceed to a definitive phase III trial, while if the posterior belief that $\delta < \delta^*$ is large enough, then development of the treatment will be abandoned. The sample size is set to ensure that the posterior belief will satisfy one of these conditions, so that an inconclusive study is avoided. Cases of binary and normally distributed responses were considered in [7], and extensions to trials yielding survival data are possible [13].

## 2. A Bayesian approach to sample size determination

Suppose that patients are randomised between a control treatment, $E_0$, and k experimental treatments, $E_1, ..., E_k$. The response $y_{ij}$ of the $i^{th}$ patient on treatment $E_j$ is distributed normally with mean $\mu_j$ and precision $\nu$ (so that the variance = $\nu^{-1}$) for $i = 1, ..., n_j$; $j = 0, 1, ..., k$. To begin with, we assume that the value of $\nu$ is known. Treatment effects are denoted by $\delta_j = \mu_j - \mu_0$, $j = 1, ..., k$. Independent normal priors for $\mu_j$ are imposed:

$$\mu_j \sim N(\mu_{0j}, (q_{0j}\nu)^{-1}), \quad j = 0, 1, ..., k.$$

Denoting the mean response of patients on $E_j$ by $\bar{y}_j$, these priors lead to independent normal posteriors for $\mu_j$:

$$\mu_j \sim N(\mu_{1j}, (q_{1j}\nu)^{-1}),$$

where

$$\mu_{1j} = (\mu_{0j}q_{0j} + n_j\bar{y}_j)/(q_{0j} + n_j)$$

and $q_{1j} = q_{0j} + n_j$ represents the total amount of information, both prior and observed, available about patients' response to treatment $E_j$, $j = 0, 1, ..., k$. Thus the posterior mean, $\mu_{1j}$, is a weighted average of the prior mean, $\mu_{0j}$, and the average of the observed responses on $E_j$, $\bar{y}_j$. The joint posterior distribution for $(\delta_1, ..., \delta_k)'$ is therefore multivariate normal:

$$\begin{pmatrix} \delta_1 \\ \delta_2 \\ \vdots \\ \delta_k \end{pmatrix} \sim N\left( \begin{pmatrix} \delta_{11} \\ \delta_{12} \\ \vdots \\ \delta_{1k} \end{pmatrix}, \begin{pmatrix} (D_{11}v)^{-1} & (q_{10}v)^{-1} & \cdots & (q_{10}v)^{-1} \\ (q_{10}v)^{-1} & (D_{12}v)^{-1} & \cdots & (q_{10}v)^{-1} \\ \vdots & \vdots & \ddots & \vdots \\ (q_{10}v)^{-1} & (q_{10}v)^{-1} & \cdots & (D_{1k}v)^{-1} \end{pmatrix} \right),$$

where

$\delta_{1j} = \mu_{1j} - \mu_{10}$ and $D_{1j} = q_{1j}q_{10}/(q_{1j} + q_{10})$, $j = 1, ..., k$.

We assume that the trial is designed to yield equal information concerning each of the experimental treatments: that is $q_{11} = ... = q_{1k}$, with common value denoted by $q_1$. Thus, sample sizes will be chosen to ensure that $q_{01} + n_1 = ... = q_{0k} + n_k$. It follows that the $D_{1j}$ will be equal to one another, and their common value is denoted by $D_1$.

Two criteria for a sufficiently large sample size will be considered. To express the first of these, denote the data collected in the trial by **y**, and posterior probabilities based on these data by $P(\bullet \mid \mathbf{y})$. Let

$$\Pi_j = P(\delta_j > 0 \mid \mathbf{y}) \quad \text{and} \quad \Gamma = P(\delta_j < \delta^* \text{ for all } j = 1,...,k \mid \mathbf{y}),$$

where $\delta^* > 0$ denotes a treatment effect of clinical importance. Hence, $\Pi_j$ is the posterior probability that treatment $E_j$ shows promise and $\Gamma$ is the posterior probability that none of the experimental treatments have a clinically important effect. The first criterion is as follows.

***Criterion 1:*** *The sample size should be sufficient to ensure that either $\Pi_j \geq \eta$ for at least one $j = 1, ..., k$, or $\Gamma \geq \zeta$, for any possible outcome dataset **y***.

If the first of the possibilities guaranteed by Criterion 1 arises, then one or more of the treatments $E_j$ for which $\Pi_j \geq \eta$ will be developed further; if the second is true, then all of the treatments will be abandoned. The values of $\eta$ and $\zeta$ will be large, and certainly greater than 0.5. There will be a small but positive probability that a dataset with the specified sample size will result in both $\Pi_j \geq \eta$ for at least one $j = 1, ..., k$ and $\Gamma \geq \zeta$ being true. In such a case, at least one treatment is showing promise but none appear to have an important effect. The sample size is adequate for the decision as to whether to study some treatment further, but further considerations such as safety or cost will be needed to reach a conclusion.

In order to calculate a suitable sample size before the study is conducted, consider a borderline outcome in which the posterior mean value for $\delta_j$ is $\delta_{1j} = \delta'_1$, $j = 1, ..., k$, leading to $\Pi_1 = ... = \Pi_k = \eta$ and $\Gamma = \zeta$. If $\delta_{1j} > \delta'_1$ for any j then $\Pi_j > \eta$ and further testing of treatment $E_j$ will be recommended. If $\delta_{1j} \leq \delta'_1$ for all j, with strict inequality for at least one, then $\Gamma > \zeta$ and no further testing of any treatment will be recommended. This situation is represented for a particular case with k = 2 in Figure 1, where the borderline outcome occurs at the single point $(\delta_{11}, \delta_{12})$ in the space of possible posterior mean treatment effects at which the *Abandon boundary* and the *Proceed boundary* meet.

In the case of general k, because $\Pi_j = 1 - P(\delta_j < 0 | \mathbf{y}) = \Phi(\delta'_1 \sqrt{D_1 v})$, it follows that $\delta'_1 = z_\eta (D_1 v)^{-1/2}$, $j = 1, ..., k$, where $z_\eta = \Phi^{-1}(\eta)$. Consider the random vector $\mathbf{X} \equiv (X_1, ..., X_k)'$ following a normal distribution with mean $\mathbf{0}$ and variance-covariance matrix $\Sigma$ for which all diagonal elements are equal to 1 and all off-diagonal elements are equal to $\rho \in [0, 1]$. Let $x_{\rho,\zeta,k}$ denote the value such that $P(\max(X_1, ..., X_k) < x_{\rho,\zeta,k}) = \zeta$. The value of $x_{\rho,\zeta,k}$ can easily be computed, or found from software such as the R program of Genz et al. [14]. Now, for this borderline situation,

$$\Gamma = P\left(\max(\delta_1, ..., \delta_k) < \delta^* | \mathbf{y}\right)$$
$$= P\left(\max\left((\delta_1 - \delta'_1)\sqrt{D_1 v}, ..., (\delta_k - \delta'_1)\sqrt{D_1 v}\right) < (\delta^* - \delta'_1)\sqrt{D_1 v} | \mathbf{y}\right) = \zeta$$

and so $(\delta^* - \delta'_1)\sqrt{D_1 v} = x_{\rho,\zeta,k}$ with $\rho = (1 + r)^{-1}$ where $r = q_{10}/q_1$. The quantity r is the ratio of posterior information required about the control treatment to that for the experimental treatments, and is set by the investigators. It follows that $\left(\delta^* - z_\eta (D_1 v)^{-1/2}\right)\sqrt{D_1 v} = x_{\rho,\zeta,k}$.

Therefore, sample sizes $n_0, n_1, ..., n_k$ should be chosen so that so that $D_{11} = ... = D_{1k} = D_1 = V_1/\nu$, where $V_1 = \{(z_\eta + x_{\rho,\zeta,k})/\delta^*\}^2$. This means that we require

$$\frac{1}{D_1} = \frac{q_1 + q_{10}}{q_1 q_{10}} = \frac{1}{q_{00} + n_0} + \frac{1}{q_{0j} + n_j} = \frac{\nu}{V_1}, \quad j = 1, ..., k. \tag{1}$$

It follows that the minimum suitable total sample size is achieved through the choices

$$n_j = (1 + k^{-1/2})V_1/\nu - q_{0j}, \, j = 1, ..., k, \text{ and } n_0 = (1 + k^{1/2})V_1/\nu - q_{00}, \tag{2}$$

and that the optimal allocation ratio r for Criterion 1 is $r = (n_0 + q_{00})/(n_j + q_{0j}) = k^{1/2}$. This is consistent with the corresponding approximate frequentist result that the optimal allocation of patients satisfies $n_0/n_j = k^{1/2}$ [15, 16].

To compute the values of $\Pi_j$ and $\Gamma$ at the end of a study, when the $q_{1j}$ may not turn out to be equal, notice that

$$\Pi_j = 1 - P(\delta_j < 0 | \mathbf{y}) = \Phi\left(\delta_{1j}\sqrt{D_{1j}\nu}\right), \tag{3}$$

and

$$\begin{aligned}
\Gamma &= P\left(\text{all of } \delta_1, ..., \delta_k \text{ are } < \delta^* | \mathbf{y}\right) \\
&= \sqrt{q_{10}\nu} \int_{-\infty}^{\infty} \left[\prod_{j=1}^{k} \Phi\left\{(\delta^* + \mu_0 - \mu_{1j})\sqrt{q_{1j}\nu}\right\}\right] \phi\left\{(\mu_0 - \mu_{10})\sqrt{q_{10}\nu}\right\} d\mu_0 \\
&= \int_{-\infty}^{\infty} \left[\prod_{j=1}^{k} \Phi\left\{u\sqrt{q_{1j}/q_{10}} + (\delta^* - \delta_{1j})\sqrt{q_{1j}\nu}\right\}\right] \phi(u) du.
\end{aligned} \tag{4}$$

Alternatively, $\Gamma$ can be calculated directly from the joint distribution of $\delta_1, ..., \delta_k$ given above, using a multivariate normal package such as [14]. In equations (3) and (4), the possibility that the $D_{1j}$ and the $q_{1j}$ might not be equal is accommodated, as sample sizes might not turn out to be quite the same as specified by the design.

A second sample size criterion involves the quantity $\Pi^*$ defined as

$$\begin{aligned}
\Pi^* &= P(\delta_j > 0 \text{ for at least one } j = 1, ..., k | \mathbf{y}) \\
&= 1 - \int_{-\infty}^{\infty} \left[\prod_{j=1}^{k} \Phi\left\{u\sqrt{q_1/q_{10}} - \delta_{1j}\sqrt{q_1\nu}\right\}\right] \phi(u) du,
\end{aligned} \tag{5}$$

and is expressed as follows.

***Criterion 2:*** *The sample size should be sufficient to ensure that either $\Pi^* \geq \eta$ or that $\Gamma \geq \zeta$ for any possible outcome dataset **y**.*

If the first of the possibilities guaranteed by Criterion 2 arises, then one or more of the treatments $E_j$ will be developed further; if the second is true, then all of the treatments will be abandoned. As for Criterion 1, the values of $\eta$ and $\zeta$ will be greater than 0.5. There will be a small but positive probability that a dataset with the specified sample size will result in both $\Pi^* \geq \eta$ and $\Gamma \geq \zeta$ being true, in which case further considerations will be needed to decide whether to proceed to further clinical study.

In order to calculate a suitable sample size, consider the limiting form of dataset **y** in which all of the $\delta_{1j}$ are equal to $-\infty$, except for $\delta_{1h}$, $\Pi^* = P(\delta_h > 0 \mid \mathbf{y}) = \Phi(\delta_{1h}\sqrt{(D_1 v)})$ and $\Gamma = P(\delta_h < \delta^* \mid \mathbf{y}) = \Phi((\delta^* - \delta_{1h})\sqrt{(D_1 v)})$. Criterion 2 is then satisfied if $\sqrt{(D_1 v)} = (z_\eta + z_\zeta)/\delta^*$. To achieve this result for all situations in which all but one of the $\delta_{1j}$ are equal to $-\infty$, we choose $n_0, n_1, ..., n_k$ so that $D_{11} = ... = D_{1k} = D_1 = V_2/v$, where $V_2 = \{(z_\eta + z_\zeta)/\delta^*\}^2$. This means that the sample sizes should be chosen to ensure that

$$\frac{1}{D_1} = \frac{q_1 + q_{10}}{q_1 q_{10}} = \frac{1}{q_{00} + n_0} + \frac{1}{q_{0j} + n_j} = \frac{v}{V_2}, \quad j = 1, ..., k. \tag{6}$$

Optimal sample sizes will satisfy (2), with $V_1$ replaced by $V_2$. This limiting case is represented for $k = 2$ in Figure 2, as the asymptotic values of $(\delta_{11}, \delta_{12})$ where the *Abandon boundary* and the *Proceed boundary* meet.

For sample sizes as specified by (6), Criterion 2 will be fulfilled in general, as will now be explained. Consider the random vector $\mathbf{X} \equiv (X_1, ..., X_k)'$ following a normal distribution, now with mean $\boldsymbol{\mu}$ that can take any value, and with variance-covariance matrix $\Sigma$ for which all diagonal elements are equal to 1 and all off-diagonal elements are equal to $\rho \in [0, 1]$. It is shown in [17] that if $P(X_1, ..., X_k < 0) \geq \kappa$ for some $\kappa \in (0, 1)$, then $P(X_1, ..., X_k <$ $z_\zeta - z_\kappa) > \zeta$ for all $\zeta \in (\kappa, 1)$. Suppose that $\Pi^* \leq \eta$. Then $P(\delta_1, ..., \delta_k < 0) \geq 1 - \eta$. The result from [17] can be applied to the posterior distribution of $\delta$, identifying $X_j$ with $\delta_j(D_1\nu)^{1/2}$, $j = 1, ..., k$. It follows that $P(\delta_1, ..., \delta_k < (z_\zeta + z_\eta)(D_1\nu)^{-1/2}) > \zeta$, as $\zeta \in (1 - \eta, 1)$ because both $\eta$ and $\zeta$ are assumed to be greater than 0.5. Equation (6) and the definition of $V_2$ show that $\delta^*(D_1\nu)^{1/2} = z_\eta + z_\zeta$, and so it follows that $\Gamma$ is greater than or equal to the latter probability, so that $\Gamma > \zeta$.

As $P(\max(X_1, ..., X_k) < x) < P(X_1 < x)$, it follows that $x_{\rho,\zeta,k} > z_\zeta$ and consequently $V_1 > V_2$. If sample sizes are set to ensure that Criterion 1 holds, then Criterion 2 will also be satisfied.

The use of Criterion 2 guarantees only that a positive trial outcome will coincide with a strong belief that at least one of the k experimental treatments is superior to control. It is possible that the final dataset will not support a strong belief of superiority for any particular one of the k experimental treatments, while at the same time not allowing conviction that none of them work. The use of Criterion 1 avoids such a dilemma, and although it leads to larger sample sizes, it is to be preferred for most purposes. For optimality, the ratio r should be set to $k^{1/2}$. Criterion 2 might be of value if a positive trial outcome would lead to further phase III development of all of the experimental treatments $E_j$ for which $P(\delta_j < \delta^* \mid \mathbf{y}) < \zeta$, perhaps concentrating on a longer-term and more definitive endpoint.

## 3. A numerical example

Suppose that k = 2 experimental treatments are to be compared with control, and that the known precision of the responses to be observed is $\nu = 1$. The prior distribution for the mean response of each of the experimental treatments is normal with mean 0.25 and precision = 4. The corresponding prior distribution for the control treatment is normal with mean 0 and precision 16. The clinically important difference is set to be $\delta^* = 0.5$. Thus, it is expected that the experimental treatments will be superior to the control by about half of the interesting difference: the prior means will be used in the Bayesian analysis of the trial although they play no part in determining the sample size. Prior knowledge of responses to the control treatment is being taken to be more reliable than that for the experimental treatments.

We set $\eta = 0.95$ and $\zeta = 0.90$. The optimal ratio of experimental to control information is $r = 2^{1/2} = 1.4142$, and this choice leads to a correlation between posterior treatment effects of $\rho = (1 + r)^{-1} = 0.4142$. The value of $x_{0.4142, 0.90, 2}$ is 1.5915. Hence $V_1 = \{(1.6449 + 1.5915)/0.5\}^2 = 41.90$. This leads to optimal sample sizes of $n_1 = n_2 = (1 + 2^{-1/2}) 41.90/1 - 4 = 67.52$ and $n_0 = (1 + 2^{1/2}) 41.90/1 - 16 = 85.15$. These values can be rounded up to $n_1 = n_2 = 68$ and $n_0 = 86$: a total sample size of 222. Figure 1 shows the decision

FIGURE 1 NEAR HERE

*(Although Figure 1 is referred to earlier in the paper, the main account of it is here.)*

boundary for the final analysis. The value of the vector $(\delta_{11}, \delta_{12})$ made up of the posterior means for the two treatment effects can be plotted on this diagram. If the plotted point lies to the right of the vertical (green) boundary, then $\Pi_1 \geq \eta$, and it can be concluded that treatment $E_1$ shows promise and further clinical study of it is appropriate. If the plotted point lies above the horizontal (green) boundary, then $\Pi_2 \geq \eta$, and treatment $E_2$ should be developed further. If both of these conditions apply, then investigators will need to choose which treatment to investigate, or whether to proceed further with both. Finally, if the plotted point is below and to the left of the curved (red) boundary, then $\Gamma \geq \zeta$, and it can be concluded that at neither treatment is likely to achieve a clinically relevant effect, and that research into them should be abandoned. The vertical green boundary is situated at $\delta'_1 = z_\eta (D_1 v)^{-1/2}$ which is equal to 0.2537 in this case. The sample sizes have been chosen to ensure that there are no points for which for which none of $\Pi_1 \geq \eta$, $\Pi_2 \geq \eta$, or $\Gamma \geq \zeta$ occur. However, there are points lying between the two boundaries for which either both $\Pi_1 \geq \eta$ and $\Gamma \geq \zeta$ or $\Pi_2 \geq \eta$ and $\Gamma \geq \zeta$ occur. If ($\delta_{11}$, $\delta_{12}$) lies in these regions, then further considerations about the treatments will have to be used to decide whether to take this line of clinical research further.

Criterion 2 depends on the value $V_2 = \{(1.645 + 1.282)/0.5\}^2 = 34.26$. From expression (5), optimal sample sizes are $n_1 = n_2 = (1 + 2^{-1/2}) \times 34.26 - 4 = 54.48$ and $n_0 = (1 + 2^{1/2}) \times 34.26 - 16 = 66.70$. Practical sample sizes are $n_0 = 67$ and $n_1 = n_2 = 55$, totalling 177. Figure 2 shows the corresponding decision boundary. The sample size that follows

FIGURE 2 NEAR HERE

*(Although Figure 2 is referred to earlier in the paper, the main account of it is here.)*

from Criterion 2 is smaller than that for Criterion 1, although as discussed earlier, this reduction in sample size comes at the cost of a risk of indecision concerning the choice of which treatment to take forward. In this context, the Bayesian procedure based on Criterion 1 does make an adjustment for the testing of multiple treatments.

For each criterion, the design recommendations involve the collection of more data on the control treatment than on each of the experimental treatments, as the control features in the estimation both treatment effects. However, the fact that more prior information is available concerning the control reduces the imbalance of randomisation. Consider Criterion 1 again, with $q_{01} = q_{02} = 4$ as before. Suppose that there is so much prior experience with the control that $q_{00} = 102$. Then the optimal samples sizes will be $n_1 = n_2 = 68$ and no observations on control ($n_0 = 0$), leading to a total sample size of 136. This situation would arise if the prior standard deviation for $\mu_0$ was 0.099, whereas those for $\mu_1$ and $\mu_2$ were 0.50. Being so much better informed about response to a standard form of treatment, relative to the uncertainty about a completely new drug would appear to be quite possible.

The need to achieve integer sample sizes has introduced approximation into the sample size calculations above, with "ideal" but non-integer solutions being rounded up for practical use. Occasionally, an exhaustive integer search might lead to a smaller total sample size than follows from rounding up an optimal solution. For example, consider again the case of Criterion 1, for which we have found a suitable design with $(n_0, n_1, n_2) = (86, 68, 68)$, and a total sample size equal to 222. In fact, five sample size combinations satisfy Criterion 1 with a total sample size equal to 221. These are $(n_0, n_1, n_2) = (85, 68, 68), (81, 70, 70), (83, 69, 69), (87, 67, 67)$ and $(89, 66, 66)$. However, it is not possible to obtain a suitable design with a total sample size less than the smallest integer exceeding $(1+k^{1/2})^2 V_1/\nu - (q_{00} + q_{01} + ... + q_{0k})$, which is the sum of the unrounded sample sizes given by (2).

Comparisons with frequentist approaches can be made. Suppose that the trial is considered as two separate comparisons, one of $E_1$ with $E_0$, and the other of $E_2$ with $E_0$. For j = 1 and j = 2, it will be concluded that $E_j$ is worthy of further research if $E_j$ is found to be significantly superior to $E_0$ at the 5% significance level (one-sided). We desire a power of 90% of drawing such a conclusion if $\delta_j = 0.5$. If the optimal allocation of m patients to each experimental treatment and $m_0 = m\sqrt{2}$ to the control is made, then the information available for each pairwise comparison will be $m_0 m/(m_0 + m) = m/(1 + 2^{-1/2}) = 0.586m$. This has to be set equal to $\{(z_{0.975} + z_{0.90})/\delta^*\}^2 = 34.27$. Hence, m = 58.50 and $m_0 = 82.74$, which would be rounded up to m = 59 and $m_0 = 83$ (a total sample size of 201). Notice that expression (6) gives precisely this solution for the Bayesian case based on Criterion 2 when there is no prior information on any of the treatments ($q_{00} = q_{01} = q_{02} = 0$). The frequentist approach here makes no adjustment for multiplicity.

Dunnett's (1984) procedure [1] does make allowance for multiplicity. The frequentist power requirements of Dunnett's method in the setting of this section can be specified as follows. If $\delta_1 = \delta_2 = 0$, then the probability of proceeding to further research with either experimental treatment should be limited to 0.05. If $\delta_1 = 0.5$ and $\delta_2 = 0$, then the probability of proceeding to further research with $E_1$ should be 0.90 (with a corresponding requirement in the case $\delta_1 = 0$ and $\delta_2 = 0.5$). For an optimal sample size, with $\sqrt{2}$ times as many patients on the experimental as on the controls, we need $n_0 = 100$ and $n_1 = n_2 = 71$, a total sample size of 242. These sample sizes are very close to the values of $(n_0 + q_{00}) = 101$ and $(n_1 + q_{10}) = 72$ that follow from Criterion 1. Even greater sample sizes result from the more general form of Dunnett's procedure in which the probability of proceeding to further research with $E_1$ should be 0.90 if $\delta_1 = 0.5$ and $\delta_2$ is equal to some non-interesting treatment effect lying between 0 and 0.5.

## 4. Allowing for uncertainty concerning the precision of patient responses

So far in this paper, it has been assumed that the precision $\nu$ of patient responses is known. In practice this will not be true. In frequentist sample size calculations, a specific value has to be assumed, in a "what if" style of argument. If the assumed value turns out to be wrong, then the intended power might not be achieved. Strategies such as sample size reviews [18, 19, 20] can be employed to reduce the risk of losing power. The Bayesian approach, on the other hand, lends itself to making proper allowance for uncertainty about the value of $\nu$.

Specification of prior opinion about $\nu$, in the form of a prior distribution, can be incorporated within the model. In particular, take this to be a conjugate gamma distribution, parameters $\alpha_0$ and $\beta_0$ leading to a posterior gamma distribution for $\nu$ with parameters $\alpha_1$ and $\beta_1$. The latter are given by $\alpha_1 = \alpha_0 + \tfrac{1}{2}n$, where n denotes the total sample size and $\beta_1 = \beta_0 + \tfrac{1}{2}H$, where

$$H = \sum_{j=0}^{k}\left(U_j + q_{0j}\mu_{0j}^2 - q_{1j}\mu_{1j}^2\right) \quad \text{and} \quad U_j = \sum_{i=1}^{n_j} y_{ij}^2, \; j = 0,1,...,k. \tag{7}$$

This is consistent with equation (6) of [7]. At the end of a study, regardless of how the sample size was determined, the values of $\Pi_j$, $\Pi^*$ and $\Gamma$ can be found by integrating the forms (1), (2) and (4) multiplied by the posterior gamma density, over $\nu$. The final Bayesian decision can be made allowing for uncertainty about the value of ν. However, if the sample sizes indicated in Section 2 for meeting Criterion 1 are used, it may now be possible to observe a dataset such that neither $\Pi_j \geq \eta$ for any j = 1, ..., k, nor $\Gamma \geq \zeta$, and correspondingly for Criterion 2. In this section, methods will be described for determining sample sizes that will satisfy the chosen criterion with high probability, rather than with certainty.

As in the known variance case, the posterior distribution of $(\delta_j - \delta_{1j})\sqrt{(D_{1j}\nu)}$, given ν, is standard normal, j = 1, ..., k. Following Appendices 2-4 of [7], it can be seen that $J_1 = 2\beta_1\nu$ has a posterior chi-squared distribution with $2\alpha_1$ degrees of freedom, independent of the latter distribution. Combining these two results it follows that, if $T_j$ is defined by

$$T_j = \frac{(\delta_j - \delta_{1j})\sqrt{D_1\nu}}{\sqrt{J_1/(2\alpha_1)}} = (\delta_j - \delta_{1j})\sqrt{\frac{D_1\alpha_1}{\beta_1}}, \tag{8}$$

then $T_j$ follows a t-distribution with $2\alpha_1$ degrees-of-freedom. The $T_j$ have pairwise correlations equal to ρ, and it can be shown that

$$P(\max(T_1,...,T_k) \leq t) = \int_0^\infty \int_{-\infty}^\infty \left\{\Phi\left(\frac{t\sqrt{v/(2\alpha_1)} + u\sqrt{\rho}}{\sqrt{1-\rho}}\right)\right\}^k \phi(u)g(v)\,du\,dv, \tag{9}$$

where $\Phi$ and $\phi$ respectively denote the distribution and density functions of a standard normal distribution, and g denotes the gamma density with parameters $\alpha_1$ and $\frac{1}{2}$. Let $t_k(2\alpha_1, \rho, \zeta)$ denote the value such that $P(\max(T_1, ..., T_k) \leq t_k(2\alpha_1, \rho, \zeta)) = \zeta$.

To find a sample size satisfying Criterion 1, we again consider a borderline outcome in which the posterior mean value for $\delta_j$ takes the common value $\delta_{1j} = \delta_1'$, j = 1, ..., k, chosen so that $\Pi_1 = ... = \Pi_k = \eta$ and $\Gamma = \zeta$. From (8) it follows that

$$\Pi_j = P(\delta_j > 0 | \mathbf{y}) = P(T_j > -\delta_1'\sqrt{D_1\alpha_1/\beta_1} | \mathbf{y}) = P(T_j < \delta_1'\sqrt{D_1\alpha_1/\beta_1} | \mathbf{y}) = \eta.$$

Hence $\delta'_1 = t(2\alpha_1, \eta)(D_1 \alpha_1/\beta_1)^{-1/2}$, j = 1, ..., k, where t(s, ξ) denotes the ξ point on the t distribution function with s degrees-of-freedom. We also require $\Gamma \geq \zeta$, where for the borderline outcome under consideration,

$$\Gamma = P\left(\max(\delta_1,...,\delta_k) < \delta^* \mid \mathbf{y}\right)$$
$$= P\left(\max(T_1,...,T_k) < (\delta^* - \delta'_1)\sqrt{\frac{D_1 \alpha_1}{\beta_1}} \mid \mathbf{y}\right).$$

Thus, we require $(\delta^* - \delta'_1)(D_1\alpha_1/\beta_1)^{1/2} \geq t_k(2\alpha_1, \rho, \zeta)$, with $D_1$ given by

$$D_1 \geq \frac{\beta_1}{\alpha_1}\left\{\frac{t(2\alpha_1,\eta) + t_k(2\alpha_1,\rho,\zeta)}{\delta^*}\right\}^2. \tag{10}$$

As $\beta_1$ is not known until the data have been observed, inequality (10) cannot be used directly as the basis of a sample size calculation. Instead, we will seek a sample size such that the probability that (10) is true exceeds some set large value ξ.

From equation (7), $\beta_1 = \beta_0 + \frac{1}{2}H$. Following Appendix 4 of [7], it can be seen that K = νH has a chi-squared distribution with n degrees of freedom. The prior distribution of $J_0 = 2\beta_0\nu$ is chi-squared with $2\alpha_0$ degrees of freedom and independent of that of K, so that $F = (K/n)/(J_0/2\alpha_0)$ follows the F-distribution with n and $2\alpha_0$ degrees of freedom and B = nF/(2α₀ + nF) the beta distribution with parameters $\frac{1}{2}$n and $\alpha_0$. As $\beta_1 = \beta_0/(1-B)$, $\beta_1$ will be less than or equal to $\beta_0/\{1 - \text{Beta}(\frac{1}{2}n, \alpha_0, \xi)\}$ with probability ξ, where Beta(a, b, ξ) denotes the 100ξ% point on the distribution function of a beta random variable with parameters a and b. Thus, if sample sizes are chosen so that

$$D_1 \geq \frac{\beta_0}{\alpha_1}\left\{\frac{1}{1-\text{Beta}(\frac{1}{2}n,\alpha_0,\xi)}\right\}\left\{\frac{t(2\alpha_1,\eta) + t_k(2\alpha_1,\rho,\zeta)}{\delta^*}\right\}^2, \tag{11}$$

then Criterion 1 will be satisfied with probability ≥ ξ.

Denoting the right-hand side of (11) by $V_n$, the minimum suitable total sample size n satisfies

$$n = \left(1+k^{1/2}\right)^2 V_n - \sum_{j=0}^{k} q_{0j}, \tag{12}$$

with individual treatment sample sizes given by $n_j = (1 + k^{-1/2})V_n - q_{0j}$, $j = 1, ..., k$, and $n_0 = (1 + k^{1/2})V_n - q_{00}$. This is the same as expression (2) with $V_1/v$ replaced by $V_n$. A search over values of n, or application of an equation solver to (12), will lead to the minimum value of n satisfying inequality (11) with the optimal allocation, $r = (n_0 + q_{00})/(n_j + q_{0j}) = k^{1/2}$, of patients to treatment. To find the smallest possible sample size, the equation is first solved for non-integer values of n, leading non-integer values for $n_0$ and the $n_j$: these are then rounded up to integer values and the final integer value of n deduced. Values of $t_k(2\alpha_1, \rho, \zeta)$ can be found from the R program of Genz et al. [14], or else (11) can be recast as

$$t_k(2\alpha_1, \rho, \zeta) \leq \delta^* \sqrt{\frac{D_1 \alpha_1}{\beta_0} \{1 - \text{Beta}(\tfrac{1}{2}n, \alpha_0, \xi)\}} - t(2\alpha_1, \eta). \tag{13}$$

The right-hand side of (13) then can be evaluated for various total sample sizes, allocated optimally, and substituted for t in equation (9) to find the smallest n for which the result is not less than $\zeta$. Both approaches were used to find and check the calculations in Section 5 below (using SAS to implement the latter method).

Following similar arguments as led to (11), it follows that suitable sample sizes to ensure that Criterion 2 is satisfied with probability $\geq \xi$ should be chosen to ensure that

$$D_1 \geq \frac{\beta_0}{\alpha_1} \left\{ \frac{1}{1 - \text{Beta}(\tfrac{1}{2}n, \alpha_0, \xi)} \right\} \left\{ \frac{t(2\alpha_1, \zeta) + t(2\alpha_1, \eta)}{\delta^*} \right\}^2. \tag{14}$$

This expression is easier to use directly than (11), as quantiles of the univariate t-distribution are widely available in statistical software.

### 5. A case study

Carlsen et al. [2] describe a clinical trial comparing four doses of bendrofluazide and placebo administered to patients suffering from arterial hypertension. The principal endpoint was the mean reduction in systolic blood pressure ($\Delta$ bp) over 10-12 weeks treatment (mm Hg). The results of the trial of bendrofluazide are given in Table I, together with some derived statistics. All four doses were successful in reducing blood pressure, with the top dose having the greatest effect.

TABLE I NEAR HERE

The authors report that sample sizes of 50 in each of the five treatments groups were found following the method of Dunnett [1]. The frequentist power requirements used in [2] were as follows. If $\delta_1 = ... = \delta_4 = 0$, then the probability of proceeding to further research with any of the experimental treatments should be limited to 0.05. If $\delta_1 = 5$ and $\delta_2 = \delta_3 = \delta_4 = 0$, then the probability of proceeding to further research with $E_1$ should be 0.90 (with corresponding requirements should any of the other treatments show such an advantage). It was assumed that the standard deviation of the reduction in blood pressure measurements ($\sigma$) was equal to 7 mm Hg. The precise result from Dunnett's method, with all five sample sizes constrained to be equal, is to use a sample size of 47 for each group, a total sample size of 235. These sample sizes are recorded in Table II, together with those from the alternative designs that are discussed in this section.

TABLE II NEAR HERE

According to the Dunnett design based on sample sizes of 47 patients, if $E_{j*}$ is the treatment associated with the highest mean endpoint, then it should be taken forward for further study if $Z_{j*} \geq c$, where $c = 2.16$ and

$$Z_j = (1/\sigma)\sqrt{n_0 n_j / (n_0 + n_j)} (\bar{y}_j - \bar{y}_0) = 0.693(\bar{y}_j - \bar{y}_0), \quad j = 1, ..., 4. \tag{15}$$

Treatment $T_4$ is associated with the highest mean endpoint, so that $Z_{j*} = Z_4 = 10.29$. As this exceeds $c = 2.16$, the fourth dose (10.0 mg/day) should be selected for further study. However, the sample sizes all exceed 47, and they differ from one another. Furthermore, the standard deviation of the observations substantially exceeds the anticipated value of 7. In Table I, the value of $Z'_j$ is computed for each dose group, where $Z'_j$ is defined in the same way as $Z_j$, but with the actual pooled estimate of standard deviation (based on the treatment group in question and the placebo group) used in place of $\sigma$, and the true sample sizes in place of those planned at the design stage. It follows that $Z'_{j*} = Z'_4 = 5.21$. Although there is no mention in [1] of how a p-value can be computed at the end of a trial, the value given by

$$p = P(Z'_{j*} \geq 5.21 \mid \delta_1 = \delta_2 = \delta_3 = \delta_4 = 0).$$

is consistent with the approach. Denoting the observed value of $Z'_{j*}$ by $z^*$ and the estimate of $\sigma$ from the responses on treatment $E_j$ by $\sigma_j$, it can be shown that,

$$p = \sum_{j=1}^{k} \int_{-\infty}^{\infty} \left\{ \prod_{i \neq j} \Phi\left( u \frac{\sigma_j}{\sigma_i} \sqrt{\frac{n_i}{n_j}} \right) \right\} \Phi\left( u \frac{\sigma_j}{\sigma_0} \sqrt{\frac{n_0}{n_j}} - z^* \frac{\sigma_j}{\sigma_0} \sqrt{\frac{n_0 + n_j}{n_j}} \right) \phi(u) \, du.$$

For these data, $p = 0.000000184$. Notice that this approach updates the values of the standard deviations by estimation from the data, but they are still treated as fixed. No allowance is made for their estimation from the data.

Now consider the Bayesian sample size calculation for this case. Suppose that $\nu = 1/7^2 = 0.0204$, and that prior information on the four doses indicates that $\mu_0 = 0$ and $\mu_1 = \mu_2 = \mu_3 = \mu_4 = 9$, and $q_{00} = 10$ and $q_{01} = q_{02} = q_{03} = q_{04} = 2$ (as shown in Table I). Take $\eta = 0.95$ and $\zeta = 0.90$. Following [2], the clinically important difference is set at $\delta^* = 5$. The optimal ratio of control to experimental information is $r = 4^{1/2} = 2$, and this choice leads to a correlation between posterior treatment effects of $\rho = (1 + 2)^{-1} = 0.3333$. The value of $x_{0.3333, 0.90, 4}$ is 1.8886. Hence $V_1 = \{(1.6449 + 1.8886)/5\}^2 = 0.4994$. This leads to optimal

sample sizes on the experimental treatments of $(1 + 4^{-1/2})\ 0.4994/0.0204 - 2 = 34.71$ and $n_0 = (1 + 4^{1/2})\ 0.4994/0.0204 - 10 = 63.41$. To a good approximation, a total of 204 patients could be recruited, with 35 patients allocated to each of the experimental treatments and 64 to control. In a similar way, suitable sample sizes satisfying Criterion 2 can be found to be 24 patients on each of the active treatments and 41 on the control: a total of 137.

Following collection of the data in Table I, a variety of Bayesian analyses can be conducted. Assuming that the common standard deviation really is equal to 7, equations (3) and (4) result in the values $\Pi_1 = ... = \Pi_4 = 1$ and $\Gamma = 0$ respectively. The value of $\Pi^*$, used with Criterion 2, is also 1. To provide values for numerical checking of computations, we also evaluate $\Gamma^{(10)} = P(\delta_j < 10 \text{ for } j=1,...,4 | \mathbf{y})$ and $\Gamma^{(15)} = P(\delta_j < 15 \text{ for } j=1,...,4 | \mathbf{y})$, to obtain $\Gamma^{(10)} = 0.000253$ and $\Gamma^{(15)} = 0.689$. This would amount to overwhelming evidence that each of the treatment effects is greater than 0, but scepticism as to whether any of them are greater than 15.

Next, recognising that each treatment group has a different standard deviation, the estimated "sd Δ bp" from Table I is used to find a separate estimate for the precision $v_j$ of each treatment $E_j$, $j = 0, ..., 5$. In equation (3), the argument of $\Phi$ becomes $\{\delta_{1j}\sqrt{(q_{10}v_0 + q_{1j}v_j)/(q_{10}v_0 q_{1j}v_j)}\}$, in the last part of equation (4), the argument of $\Phi$ becomes $\{u\sqrt{(q_{1j}v_j)/(q_{10}v_0)} + (\delta^* - \delta_{1j})\sqrt{q_{1j}v_j}\}$, and similar changes are made to the expression for $\Pi^*$. We still obtain $\Pi_1 = ... = \Pi_4 = 1$, $\Pi^* = 1$ and $\Gamma = 0$ to four decimal places, but now $\Gamma^{(10)} = 0.0168$ and $\Gamma^{(15)} = 0.562$.

Now suppose that the uncertainty in the prediction of $v$ was anticipated. The prior mean of $v$ is taken to be $7^{-2} = 0.0204$, consistently with the value of common standard deviation used in [2]. The largest observed value of standard deviation is 15, corresponding to a precision of 0.00444. Choosing a gamma prior for $\nu$ with parameters $\alpha_0 = 1$ and $\beta_0 = 49$, provides a prior mean value of 0.0204 with a prior probability that $\nu \leq 0.00444$ of 0.196. With hindsight, this would appear to be a sensibly vague choice. Using the expressions in (7) reveals that the posterior distribution of $\nu$ is gamma with parameters $\alpha_1 = 1 + 257/2 = 129.5$ and $\beta_1 = 49 + \frac{1}{2}H$, where $H = 46413.54$ so that $\beta_1 = 23255.77$. The treatment group contributions to H are listed in Table I as $h_j$. To compute the $h_j$, the sums of squares of endpoints within each treatment group are needed: these are displayed in Table I as ss. They are deduced from the values of standard errors which are reported in [2] to only one decimal place. Consequently the sums of squares have been reproduced only approximately. The posterior mean of $\nu$ is 0.00557, corresponding to an estimated standard deviation for the observations of $\Delta$ bp of 13.40. The posterior probability that $\nu \leq 0.00444$ is 0.00729, indicating that there is little uncertainty remaining in the posterior opinion about $\nu$. Expected values of expressions (3) and (4) can be found, treating them as functions of $\nu$, and using its posterior gamma density. The resulting values are $\Pi_1 = ... = \Pi_4 = 1$, $\Pi^* = 1$ and $\Gamma = 0$ to four decimal places. It also follows that $\Gamma^{(10)} = 0.0197$ and $\Gamma^{(15)} = 0.563$. In this case, the posterior opinion about $\nu$ is so accurate that allowing for uncertainty makes little difference. A more important distinction between these calculations and the previous set is that these assume a common value of $\nu$ across treatment groups, with uncertainty about its value, whereas the earlier calculations took different values of $\nu$ for each treatment comparison, but neglected any uncertainty about their values.

The highest dose, 10 mg/day, is associated with the largest mean response, and so it should be the one to be taken forward for further research. The posterior probability that each of the other doses actually has a larger mean effect than the top dose can be found, with the individual estimated treatment group precisions being used as the $\nu_j$, but without allowance for the uncertainty of the precision: for the sample size collected this appears to be unnecessary. For doses 1.25, 2.5 and 5, the probabilities that the corresponding treatment effect is larger than that of 10 mg/day are 0.073, 0.158 and 0.112 respectively. Unless there are safety concerns or cost considerations, the focus of further research should be on the top dose. If there are such concerns, indicating that the top dose should be avoided, then Criterion 1 indicates which treatments could be admitted to further research. In this example, as each of $\Pi_2$, $\Pi_3$ and $\Pi_4 > 0.95$ (in fact they are all equal to 1), any of $T_2$, $T_3$ and $T_4$ could be developed further.

TABLE III NEAR HERE

Finally, suppose that the uncertainty in the prediction of $\nu$ were allowed for in the sample size calculation. As discussed above, a gamma prior for $\nu$ with parameters $\alpha_0 = 1$ and $\beta_0 = 49$, is adopted. Setting $\xi = 0.95$ and choosing Criterion 1 leads to the optimal sample sizes of 714 on each experimental treatment and $n_0 = 1422$ on control: a total of $n = 4278$ patients. There is a very high price to pay for certainty of reaching a convincing conclusion in the face of uncertainty about the variability of patient responses. Because of the mismatch between the anticipated standard deviation of the mean reduction in systolic blood pressure measurements of 7, and the observed values the largest of which is 15, a broad prior distribution was set for $\nu$. In Table III, sample sizes for values of $\xi$ less than the stringent setting of 0.95, and for tighter prior distributions for $\nu$ are presented, for both criteria. Bear in mind that the actual study recruited 257 patients, and so all but one of the designs shown would have called for additional subjects.

6. **Discussion**

Bayesian methods of design and analysis are attractive for early phase clinical trials in which the objective is to decide whether and in what form to take a treatment forward for further clinical research, rather than to reach a definitive conclusion on efficacy. Bayesian methods for sample size calculation complement such analyses, and this paper demonstrates that it is relatively easy to extend such calculations to the multiple treatment setting. Caution about using the approach of this paper was expressed by Zaslavsky [21] and discussed by Zaslavsky and Whitehead [22] in the context of comparisons between a single experimental treatment and a control. In the multiple treatment setting, even more care is needed in making comparisons between sample sizes calculated from frequentist and Bayesian approaches. This is because the criteria that each is set to satisfy are so different. In particular, the Bayesian approaches allow for prior opinion effectively to stand in for real prospective observations. These concerns should be less acute during the early phase of drug development: the Bayesian approach would appear to be less suited for either the design or the analysis of phase III confirmatory studies.

Whitehead et al. [7] included discussion of Bayesian sample determination when comparing a single experimental treatment with control in terms of non-normally distributed endpoints such as binary or ordinal responses. It was supposed that the advantage of the experimental treatment over the control was expressed in terms of the single scalar parameter $\theta$, and that the analysis was to be conducted in terms of the efficient score $Z$ and Fisher's information $V$. Approximately, when the amount of information is large and the value of $\theta$ is small, $Z$ is normally distributed with mean $\theta V$ and variance $V$. In the Bayesian context, a convenient conjugate prior for $\theta$ is normal with mean $Z_0/V_0$ and variance $1/V_0$, leading to a posterior that is normal with mean $(Z_0 + Z)/(V_0 + V)$ and variance $1/(V_0 + V)$. When there are multiple experimental treatments, $E_1, ..., E_k$, then the advantage of $E_j$ over control can be parameterised as $\theta_j$, and the corresponding efficient score statistic and Fisher's information denoted by $Z_j$ and $V_j$ respectively, $j = 1, ..., k$. The statistics $Z_h$ and $Z_j$ will be correlated, $h, j = 1, ..., k$, and a form for their covariance can be deduced in each special case (see [23] for the cases of binary and ordinal data). Using such results and starting with independent priors for the $\theta_j$, the method of this paper can be adapted for large trials with non-normal responses.

In this paper, two different criteria for specifying a suitable sample size have been suggested. Criterion 1 is likely to be appropriate for most purposes as it ensures that, provided the final evidence does not lead to strong evidence that none of the new treatments is superior to control, there will be at least one treatment which is strongly supported by the evidence. Criterion 2 may be of use in circumstances where it is not necessary to identify a clearly winning treatment before concluding that a longer-term definitive trial is justified. It is interesting to observe that Criterion 1 does adjust the sample size required for the multiplicity of experimental treatments, even in a Bayesian context.

Due to the Bayesian paradigm, no adjustment is required for taking multiple looks at data emerging from the trial. Thus, interim analyses can be performed to assess whether the criterion $\Gamma \geq \zeta$ is yet satisfied, and if so to stop the study for futility. The criteria $\Pi_1 \geq \eta$, ..., $\Pi_k \geq \eta$ can also be checked to see whether a winning treatment can be identified early, or if Criterion 2 is being used a further trial can be initiated as soon as it is found that $\Pi^* \geq \eta$. In fact, achievement of Criterion 2 could be used as the trigger for investing in the often lengthy preparations for a phase III trial, while the current study continues until Criterion 1 is satisfied and a winning treatment identified for that subsequent trial. Appropriate criteria for dropping individual experimental treatments early could also be set. No alteration to the maximum sample sizes per treatment will be needed in order to maintain the Bayesian properties of the overall procedure.

**Acknowledgement**

The authors are grateful to Simon Day of Clinical Trials Consulting & Training Limited and Nelson Kinnersley of Roche Products Ltd for discussion of motivational examples and practical insights that helped in the preparation of this paper.

**Table I: Data from the study of Carlsen et al., with associated statistics and Bayesian quantities**

| Dose (mg/day) | 0 | 1.25 | 2.5 | 5.0 | 10.0 |
|---|---|---|---|---|---|
| n | 52 | 50 | 52 | 52 | 51 |
| mean $\Delta$ bp | 2.8 | 12.7 | 14.3 | 13.4 | 17.0 |
| se $\Delta$ bp | 1.7 | 2.0 | 1.6 | 2.0 | 2.1 |
| sd $\Delta$ bp | 12.3 | 14.1 | 11.5 | 14.4 | 15.0 |
| ss | 8072 | 17865 | 17423 | 19945 | 25985 |
| pooled S | - | 13.2 | 11.9 | 13.4 | 13.7 |
| $Z_j$ | - | 7.14 | 8.38 | 7.72 | 10.29 |
| $Z'_j$ | - | 3.78 | 4.93 | 4.04 | 5.27 |
| $q_{0j}$ | 10 | 2 | 2 | 2 | 2 |
| $\mu_{0j}$ | 0 | 9 | 9 | 9 | 9 |
| $q_{1j}$ | 62 | 52 | 54 | 54 | 53 |
| $\mu_{1j}$ | 2.35 | 12.56 | 14.10 | 13.24 | 16.70 |
| $\delta_{1j}$ | - | 10.21 | 11.76 | 10.89 | 14.35 |
| $h_j$ | 7730 | 9826 | 6843 | 10645 | 11369 |

**Table II: Frequentist and Bayesian designs for the case of known precision $\nu = 0.0204$**
**Frequentist settings: type I error = 0.05; power = 0.90; $\delta^* = 5$**
**Bayesian settings: $q_{00} = 10$; $q_{01} = q_{02} = q_{03} = q_{04} = 2$; $\eta = 0.95$; $\zeta = 0.90$; $\delta^* = 5$ and $r = 2$**

| Design type | sample size: | | |
|---|---|---|---|
| | on experimental treatments | on control treatment | total |
| Trial conducted by Carlsen et al [2] | 50-52 | 50 | 257 |
| Frequentist Dunnett design with equal sample sizes | 47 | 47 | 235 |
| Bayesian design, Criterion 1 | 35 | 64 | 204 |
| Bayesian design, Criterion 2 | 24 | 41 | 137 |

**Table III: Bayesian designs for the case of unknown precision**
**Settings are $q_{00} = 10$; $q_{01} = q_{02} = q_{03} = q_{04} = 2$; $\eta = 0.95$; $\zeta = 0.90$; $\delta^* = 5$ and $r = 2$**

| $\alpha_0$ | $\beta_0$ | prior: | | $\xi$ | sample size (Criterion 1/Criterion 2): | | |
|---|---|---|---|---|---|---|---|
| | | expectation of $\nu$ | probability that sd $\geq 15$ | | on experimental treatments | on control treatment | total |
| 1 | 49 | 0.0204 | 0.196 | 0.95 | 714/489 | 1422/972 | 4278/2928 |
| | | | | 0.80 | 163/111 | 320/216 | 972/660 |
| | | | | 0.50 | 52/35 | 97/63 | 305/203 |
| 2 | 98 | 0.0204 | 0.071 | 0.95 | 205/140 | 403/274 | 1223/834 |
| | | | | 0.80 | 88/59 | 169/112 | 521/348 |
| 3 | 147 | 0.0204 | 0.029 | 0.95 | 133/91 | 259/175 | 791/539 |
| | | | | 0.80 | 70/48 | 134/89 | 414/281 |

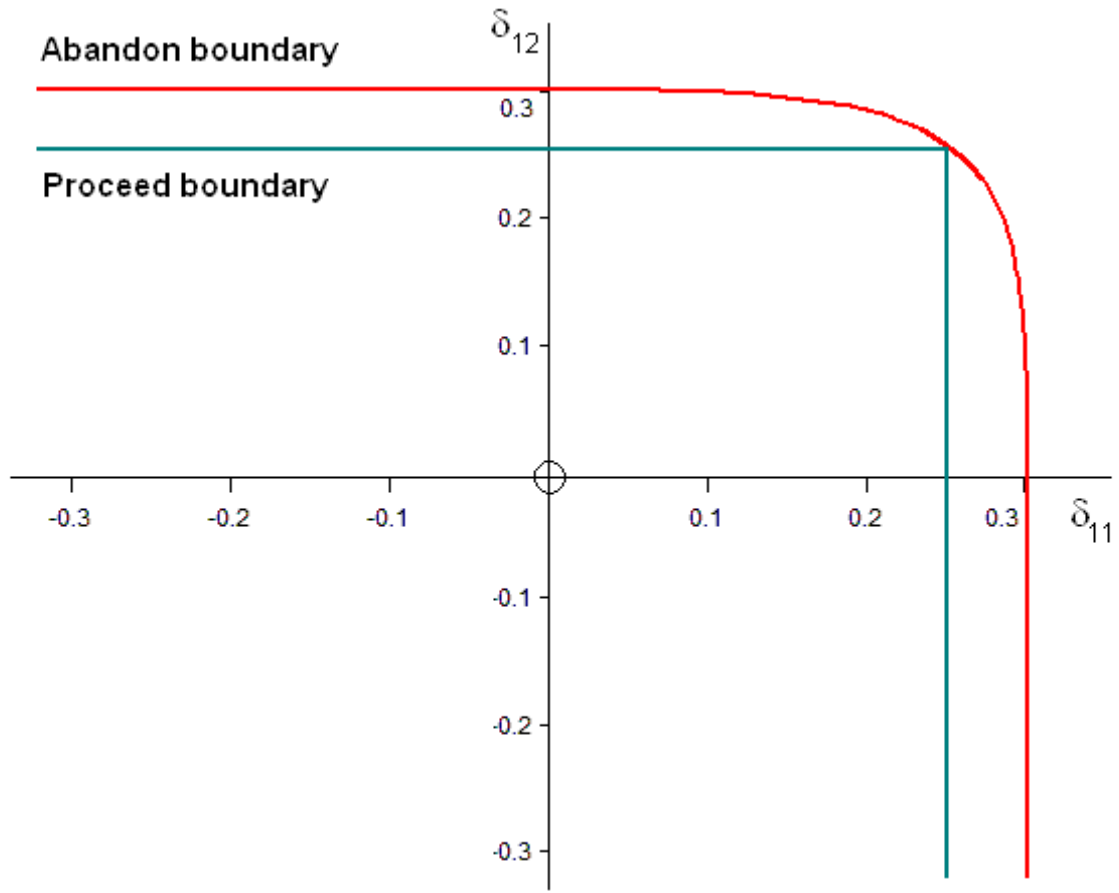

*Figure 1: An illustrative outcome space based on Criterion 1. The decision regions are: proceed to further research if ($\delta_{11}$, $\delta_{12}$) lies above or to the right of the inner "Proceed boundary"; abandon these experimental treatments if ($\delta_{11}$, $\delta_{12}$) lies below and to the left of the outer "Abandon boundary".*

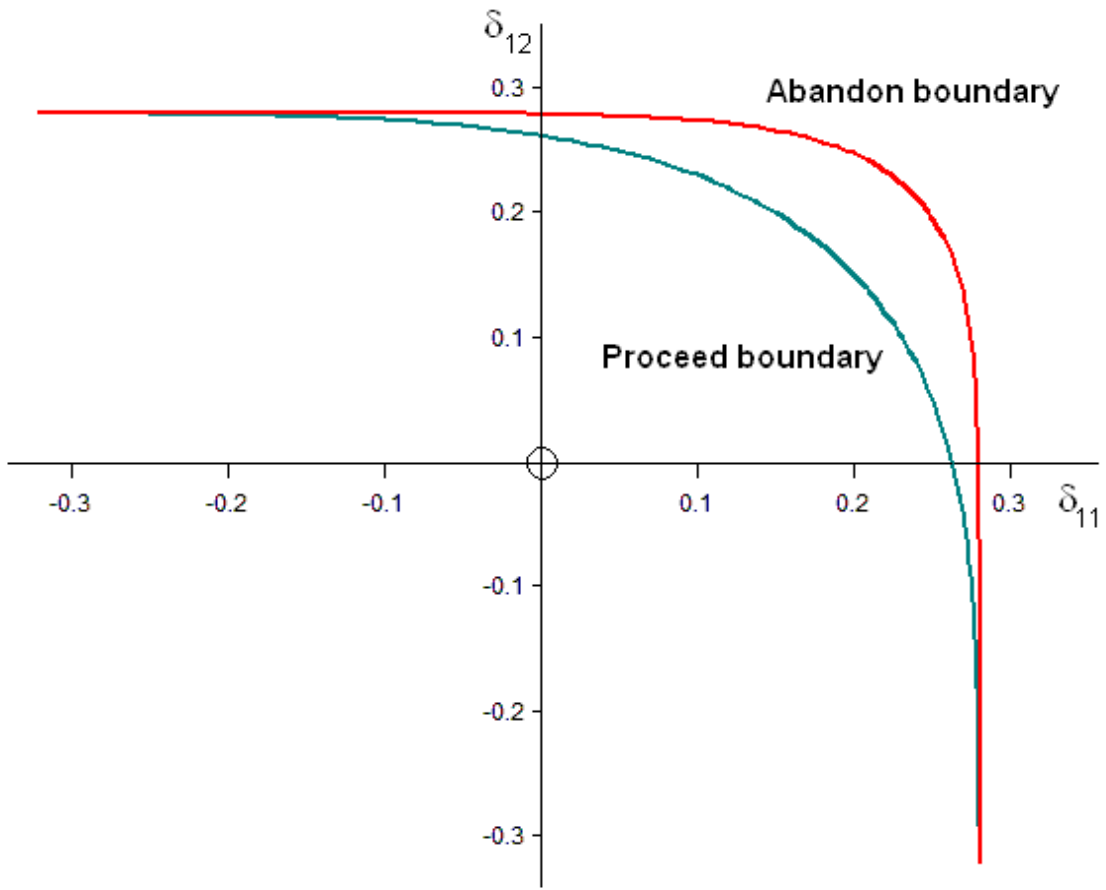

*Figure 2: An illustrative outcome space based on Criterion 2. The decision regions are: proceed to further research if ($\delta_{11}$, $\delta_{12}$) lies above or to the right of the inner "Proceed boundary"; abandon these experimental treatments if ($\delta_{11}$, $\delta_{12}$) lies below and to the left of the outer "Abandon boundary".*